\documentclass[a4paper,12pt]{amsart}

\usepackage{amssymb}
\usepackage{xy}
\xyoption{all}
\xyoption{poly}

\newtheorem{theo}{Theorem}[section]
\newtheorem{coro}[theo]{Corollary}

\newtheorem{lemma}[theo]{Lemma}
\newtheorem{defi}[theo]{Definition}
\newtheorem{conj}[theo]{Conjecture}

\def\BC{{\mathbf{C}}}

\def\BS{{\mathbf{S}}}

\def\BS{{\mathbf{S}}}

\def\BC{{\mathbb{C}}}

\def\BZ{{\mathbb{Z}}}

\def\bs{{\mathbf{s}}}

\def\Int{\operatorname{Int}\nolimits}
\def\Red{\operatorname{Red}\nolimits}
\def\Aut{\operatorname{Aut}\nolimits}

\def\ie{{\em i.e.}}

\title{A dual braid monoid for the free group}
\author{David Bessis}
\address{DMA, \'Ecole normale sup\'erieure, 
45 rue d'Ulm, 75230 Paris cedex 05, France}
\email{david.bessis@ens.fr}
\thanks{This article is the fruit of an
inspiring visit to KIAS (Seoul) in June 2003.
I thank Sang Jin Lee
for his hospitality, for stimulating discussions and for important 
suggestions.}

\begin{document}

\begin{abstract}
We construct a quasi-Garside monoid structure for the free group.
This monoid should be thought of as a \emph{dual braid monoid} for
the free group, generalising the constructions by Birman-Ko-Lee
and by the author of new Garside monoids for Artin groups of spherical type.
Conjecturally, an analog construction should be available for arbitrary
Artin groups and for braid groups of well-generated complex reflection
groups.
\end{abstract}

\maketitle

This article continues the exploration of the theory of Artin groups and
generalised 
braid groups from the new point of view introduced by Birman-Ko-Lee in
\cite{BKL} for the classical braid group on $n$ strings.
In \cite{dualmonoid}, we generalised their construction to Artin groups
of spherical type. In the current article, we study the case of the free
group, which is the Artin group associated with the universal Coxeter group.
The formal analogs of the main statements in \cite{dualmonoid} turn out
to be elementary consequences of classical material (some of which
was known to Hurwitz and Artin).
In an attempt to interpolate some recent generalisations of the dual monoid
construction (by Digne for the Artin group of type $\widetilde{A_n}$, 
\cite{digne};
by Corran and the author for the braid group of the complex reflection
group $G(e,e,n)$, \cite{BC}), we propose two conjectures describing
properties of a generalised dual braid monoid, in the contexts of
\begin{itemize}
\item[(a)] arbitrary Artin groups and
\item[(b)] braid groups of well-generated finite complex reflection groups.
\end{itemize}
This would provide the first uniform combinatorial
approach to these objects.
The initial motivation for the current work was to understand
the situation (b) from a natural geometric viewpoint;
the conjectures about complex reflection groups will
be studied in the sequel
\cite{combi}, answering some questions raised in \cite{bmr}.

\section{Hurwitz action}

For any positive integer $n$, the ``usual'' braid group is the
abstractly presented group
$$B_n:=\left< \sigma_1,\dots,\sigma_{n-1} \left| \sigma_i \sigma_{i+1}
\sigma_{i} =\sigma_{i+1}\sigma_i \sigma_{i+1}, \; \sigma_i\sigma_j
= \sigma_j\sigma_i \; \text{if} \;  |i-j| >1 \right. \right>.$$
In the problems we are interested in, two ``braid groups''
simultaneously come into play: this ``usual''
braid group, and the Artin group associated
with a Coxeter system (or the generalised braid group associated with a
complex reflection group). Except in the final conjectures, this
Artin group will be the free group.

Let $G$ be a group. For any sequence $(g_1,\dots,g_n)\in G^n$,
set $$\sigma_i \cdot (g_1,\dots,g_n) := (g_1,\dots,g_{i-1},g_ig_{i+1}g_i^{-1},
g_i,g_{i+2},\dots,g_n).$$ It is straightforward (and well-known) that
this assignment extends to a left-action of $B_n$ on $G^n$.

\begin{defi}
This action is called \emph{Hurwitz action} of $B_n$ on $G^n$.
\end{defi}

This action can be viewed as a particular example of a more general 
construction, where the important property of $G^n$ is that it is 
an \emph{automorphic set} (in the sense of \cite{brieskorn})
or equivalently a \emph{rack} (in the sense, for example, of \cite{panobraid}).

In \cite{brieskorn}, Brieskorn considers several problems about braid group
actions on automorphic sets. One of these problems is to characterise orbits.
A very naive invariant of Hurwitz action is the product
$$\pi:(g_1,\dots,g_n) \mapsto g_1\dots g_n.$$
We will be interested in situations where $\pi^{-1}(g)$ is a single
Hurwitz orbit, for a specific $g\in G$.

\section{Non-crossing loops}

In all this section, we fix $n+1$ distinct points $x_0,\dots,x_n$ in
$\BC$. The complex line is endowed with an orientation called ``positive''
or ``direct''.

We set $$F_n:=\pi_1(\BC-\{x_1,\dots,x_n\},x_0).$$
This group is isomorphic to a (``the'')
free group on $n$ generators, but its geometric definition gives
additional structure, which is what matters here.
For example, we may consider the following natural elements in $F_n$:

\begin{defi}
A \emph{non-crossing loop} is a continuous embedding
$\lambda:S^1 \hookrightarrow \BC-\{x_1,\dots,x_n\}$ whose image contains
$x_0$.

To any non-crossing loop $\lambda$, we associate the element
$f_{\lambda}\in F_n$ obtained by following $\lambda$ with the positive
orientation (coming from the orientation of $\BC$).
Elements $f_{\lambda}\in F_n$ which may be obtained this are said to 
be \emph{non-crossing}.
We denote by $NC$
the set of non-crossing elements in $F_n$
\end{defi}

We consider the length function
\begin{eqnarray*}
l:F_n & \longrightarrow & \BZ \\
f & \longmapsto & \frac{1}{2i\pi}\sum_{j=1}^n\int_f \frac{dz}{z-x_j}
\end{eqnarray*}

For any non-crossing loop $\lambda$, we may consider
the set $\Int(\lambda)$ of points of $\BC$ which are ``inside''
$\lambda$ (in the weak sense: we consider the support of $\lambda$
to be ``inside'').
Clearly, the index of $f_{\lambda}$ around $x_i$ is
$1$ if $x_i\in \Int(\lambda)$, $0$ otherwise.
Setting $ht(\lambda) := | \Int(\lambda) \cap \{x_1,\dots,x_n\}|$, we have
the relation
$$l(f_{\lambda}) = ht(\lambda).$$

\begin{defi}
We define a relation $\subseteq$ in $NC$ by
$$\forall f, g\in NC, f \subseteq g
\stackrel{\text{def}}{\Longleftrightarrow}
\exists \; \text{non-crossing loops} \;\lambda,\mu, f=f_{\lambda}, g = f_{\mu},
\Int(\lambda) \subseteq \Int(\mu).$$
\end{defi}

We leave to the reader the following easy topological lemma:

\begin{lemma}
\label{stupidlemma}
For all $f,g\in NC$, the following assertions are equivalent:
\begin{itemize}
\item[(i)] $f\subseteq g$;
\item[(ii)] for any non-crossing loop $\lambda$ such that
$f=f_{\lambda}$, there exists a non-crossing loop $\mu$ such that
$g=f_{\mu}$ and $\Int(\lambda) \subseteq \Int(\mu)$;
\item[(iii)] for any non-crossing loop $\mu$ such that
$g=f_{\mu}$, there exists a non-crossing loop $\lambda$ such that
$f=f_{\lambda}$ and $\Int(\lambda) \subseteq \Int(\mu)$.
\end{itemize}
\end{lemma}

\begin{lemma}
\begin{itemize}
\item[(i)] For all $f,g\in NC$, $f\subseteq g$ implies $l(f) \leq l(g)$.
If $f\subseteq g$ and $l(f) = l(g)$, then $f=g$.
\item[(ii)]
The relation $\subseteq$ is an order relation.
\end{itemize}
\end{lemma}

\begin{proof}
(i): The first statement is trivial. For the second statement,
choose $\lambda$ and $\mu$ such that
$\lambda,\mu, f=f_{\lambda}, g = f_{\mu}$. Since $ht(\lambda)=
ht(\mu)$, the annulus ``between'' $\lambda$ and $\mu$ contains
no point in $\{x_1,\dots,x_n\}$, thus $\lambda$ and $\mu$ are
isotopic.

(ii):
The relation is clearly reflexive. Antisymmetry follows from (i).
Transitivity follows from Lemma \ref{stupidlemma}.
\end{proof}

The main result of this section says that certain subposets of
$NC$ are lattices. Before stating it, let us observe that
$NC$ as a whole
is not a lattice. A first obstruction is that one may find
non-isotopic height $n$ non-crossing loops. Clearly, they do not
even have a common upper bound (let alone a least common upper bound).
For $n=2$,
two such loops are illustrated below (one with a full line, the other one
with a dotted line):

$$
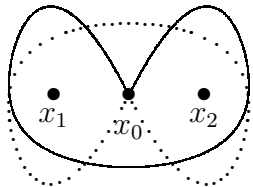
\xy
(0,-4.5) *++{x_0}, (-10,-3) *++{x_1}, (10,-3) *++{x_2},
(-10,0) *++{\bullet}, (10,0) *++{\bullet}, (0,0) *++{\bullet},
(0,0)="0", (-10,0)="1", (20,0)="2",
"0";"0" **\crv{(-10,20)&(-20,0)&(-10,-10)&(10,-10)&(20,0)&(10,20)}, 
"0";"0" **\crv{~*=<4pt>{.}(-10,-20)&(-20,0)&(-10,10)&(10,10)&(20,0)&(10,-20)}
\endxy$$

We may also observe that the corresponding elements in $NC$ do not have
a largest common lower bound:
the two height $1$ non-crossing loops represented below are distinct
maximal common
lower bounds to the above height $2$ non-crossing loops:

$$
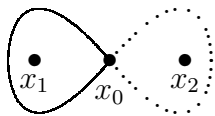
\xy
(0,-4.5) *++{x_0}, (-10,-3) *++{x_1}, (10,-3) *++{x_2},
(-10,0) *++{\bullet}, (10,0) *++{\bullet}, (0,0) *++{\bullet},
(0,0)="0", (-10,0)="1", (20,0)="2",
"0";"0" **\crv{(-10,12)&(-17,0)&(-10,-12)}, 
"0";"0" **\crv{~*=<4pt>{.}(10,12)&(17,0)&(10,-12)}
\endxy$$

\begin{defi}
For any $g\in NC$, we set $NC_g:=\{f\in NC | f \subseteq g\}$.
\end{defi}

\begin{theo}
\label{theolattice}
For any $g\in NC$, the poset $(NC_g,\subseteq)$ is a lattice.
\end{theo}

The author thanks Sang Jin Lee, for suggesting to use hyperbolic geometry in
the following proof.

\begin{proof}
First, it is easy to reduce the question to the case when $l(g)=n$.

Up to isotopy, we may assume that $x_0=-1$, and that $g$ is represented
by the unit circle. We set 
$$D:= \{ z\in \BC \mid |z| \leq 1 \}.$$
Using Lemma \ref{stupidlemma} (ii) and (iii), we observe that we may forget
the outside of $D$: any element of $f\in NC_g$ is represented by 
non-crossing loops $\lambda$ with
$\Int(\lambda)\subseteq D$, and in $NC_g$ the relation $\subseteq$ could
be equivalently redefined using only such loops.

If $n=1$, the result is straightforward.

Assume now that  $n>1$. We may endow
$D_n:=D-\{x_1,\dots,x_n\}$ with a complete hyperbolic
metric (see, for example, \cite{panobraid}, Chapter 7)
Let $\widetilde{D_n}$ be the universal cover of $D_n$ may be viewed
as a subset of the hyperbolic plane (see the nice picture on page 114,
\emph{loc. cit.}).

Any element $f\in F_n$ may be represented by a (possibly self-intersecting)
loop
in the pointed space $(D_n,x_0)$, thus be a path in $\widetilde{D_n}$;
among such paths, there is a unique geodesic. The corresponding
loop in $(D_n,x_0)$ is called the \emph{geodesic loop} of $f$.
Geodesic loops minimise self-intersections and mutual intersections;
in particular:
\begin{itemize}
\item For all $f\in F_n$, then $f\in NC_g$ if and only if its
geodesic loop is non-crossing.
\item For all $f,f'\in NC_g$ with geodesic loops $\lambda,\lambda'$,
then $f\subseteq f' \Leftrightarrow \Int(\lambda) \subseteq \Int(\lambda')$.
\end{itemize}
The theorem is a trivial consequence of the last statement:
Let $f,f'\in NC_g$ with geodesic loops $\lambda,\lambda'$.
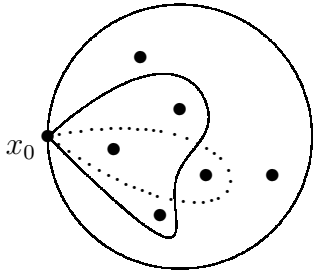
$$\xy
0;<5pt,0pt>:
(-12,-1) *++{x_0},
(-10,0) *++{\bullet},(-5,-1) *++{\bullet},
(-1.5,-6) *++{\bullet}, (2,-3) *++{\bullet}, (0,2) *++{\bullet},
(-3,6) *++{\bullet},(7,-3) *++{\bullet},
(-10,0)="0", (-10,0)="1", (20,0)="2",
"0";*\xycircle(10,10){},
"0";"0" **\crv{(2,-12)&(-3,-2)&(5,0)&(-1,8)}, 
"0";"0" **\crv{~*=<4pt>{.}(8,-13)&(9,4)}
\endxy$$
Any $h\in NC_g$ such that $f\subseteq h$ and $f'\subseteq h$ may
be represented by a non-crossing loop $\nu$ such that
$\Int(\lambda)\subseteq \Int(\nu)$ and $\Int(\lambda')\subseteq \Int(\nu)$.
Consider the loop $\lambda \vee \lambda'$ obtained by glueing the successive
``outermost'' portions of the two loops (in the above example, this
element is made with three successive portions of loops).
Clearly, any non-crossing loop containing $\Int(\lambda)\cup \Int(\lambda')$
in its interior must also contain $\lambda \vee \lambda'$ in its interior:
the element represented by $\lambda \vee \lambda'$ is the minimal
least upper bound of $f$ and $f'$.

Similarly, considering the connected component of
$\Int(\lambda)\cap \Int(\lambda')$ containing $x_0$, we obtain
a maximum lower bound. An illustration with the above $f,f'$ is given
below (the original loops are the dotted curves, the $\inf$ and the
$\sup$ are the full curves).
$$
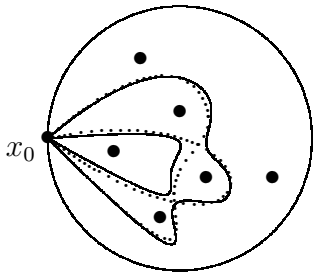
\xy
0;<5pt,0pt>:
(-12,-1) *++{x_0},
(-10,0) *++{\bullet},(-5,-1) *++{\bullet},
(-1.5,-6) *++{\bullet}, (2,-3) *++{\bullet}, (0,2) *++{\bullet},
(-3,6) *++{\bullet},(7,-3) *++{\bullet},
(-10,0)="0", (-10,0)="1", (20,0)="2",
"0";*\xycircle(10,10){},
"0";"0" **\crv{~*=<3pt>{.}(2,-12)&(-3,-2)&(5,0)&(-1,8)}, 
"0";"0" **\crv{~*=<3pt>{.}(8,-13)&(9,4)},
"0";"0" **\crv{(2,-12)&(-3,-3)&(4,-6)&(4,-2)&(1,-1)&(4,4)&(-3,6)},
"0";"0" **\crv{(1.4,-6.4)&(-3.5,-2)&(3.1,1)}
\endxy$$

\end{proof}

{\bf \flushleft Remark.} In the last proof, instead of using hyperbolic
geometry, one could use a more computational viewpoint, which may also be
used to implement the $\inf$ and $\sup$ operations.
Say that two non-crossing loops are \emph{tight} if their number of 
intersections is minimal (within their homotopy classes).
A first observation is that \emph{tight} representatives for a pair
of elements of $NC_g$ may be obtained by successive ``bigon eliminations'':
a \emph{bigon} is portion of the picture looking like
$$
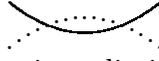
\xy
(-10,3)="0", (-10,-3)="1", (10,-3)="2", (10,3)="3",
"0";"3" **\crv{(0,-5)}, 
"1";"2" **\crv{~*=<4pt>{.}(0,5)}
\endxy$$
with no marked point $x_i$ in the inside portion; eliminating such
a bigon consists of replacing this portion of the picture by something like
$$
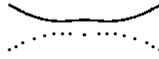
\xy
(-10,3)="0", (-10,-3)="1", (10,-3)="2", (10,3)="3",
"0";"3" **\crv{(-5,0)&(0,2)&(5,0)}, 
"1";"2" **\crv{~*=<4pt>{.}(-5,0)&(0,-2)&(5,0)}
\endxy$$
Tightness may be detected by the absence of bigons.
One may actually prove (by bigon elimination)
the stronger result: for any triple of non-crossing loops,
one may find homotopic loops which are pairwise tight.
The only property of hyperbolic geodesics used above is that they
are pairwise tight, thus that they solve the latter problem.
However, for practical use, it is very efficient to perform 
bigon elimination without relying on hyperbolic geometry.

\section{Braid reflections and coordinate systems}

Since $F_n$ is
the fundamental group of the complement in $\BC$ of a complex algebraic
hypersurface (a finite set), we may consider special elements usually
called \emph{generators-of-the-monodromy} or \emph{meridiens}
(we prefer here to call them \emph{braid reflections}).

These elements may be described as follows.
A \emph{connecting path} is a continuous map
$\gamma:[0,1] \rightarrow \BC$ such that
$\gamma(0)=x_0$,
$\gamma(1) \in \{x_1,\dots,x_n\}$ and
$t\neq 1 \Rightarrow \gamma(t) \notin \{x_1,\dots,x_n\}$.
One may associate to such a $\gamma$ an element $r_{\gamma}$ as follows:
starting from $x_0$, follow $\gamma$; arriving close to $\gamma(1)$,
make a positive turn around a small circle centered on $\gamma(1)$;
return to $x_0$ following $\gamma$ backwards.

\begin{defi}
An element $r\in F_n$ is a \emph{braid reflection} if there exists a
connecting path $\gamma$
such that $r=r_{\gamma}$.
The set of reflections in $F_n$ is denoted by $R$.
\end{defi}

\begin{lemma}
The set $R$ coincides with the set of non-crossing elements of 
height $1$.
\end{lemma}

\begin{proof}
If $r$ is non-crossing of height $1$, then choose a non-crossing loop
$\lambda$ representing $r$. We have
$\Int(\lambda)\cap \{x_1,\dots,x_n\} = \{x_{i_0}\}$. 
Since $\Int(\lambda)$ is path connected, we may draw inside $\lambda$
a path $\gamma$
connecting $x_0$ and $x_{i_0}$. It is clear that $r =r_{\gamma}$.

To prove the converse statement, one may check that for 
any path $\gamma$ connecting $x_0$ and some $x_i$, there exists
$\tilde{\gamma}$ without self-intersections such
that $r_{\gamma}=r_{\tilde{\gamma}}$ (it is clear by construction
that $r_{\tilde{\gamma}}$ is non-crossing of height $1$).
To find such a $\tilde{\gamma}$,
one may remove self-intersections by ``sliding'' them past $x_0$.
[Alternatively, one could observe that the conjugacy classes in $R$
are indexed by the irreducible components of the hypersurface; that
each conjugacy class contains a non-crossing element; and finally that
$NC$ is stable under conjugacy.]
\end{proof}

The standard way to see $F_n$ as an abstractly presented group
(with $n$ generators and no relation) is by means of a \emph{coordinate
system}:

\begin{defi}
Consider a planar graph $\Gamma$, whose vertices are $x_0,\dots,x_n$, 
and with $n$ edges $\gamma_1,\dots,\gamma_n$, each $\gamma_i$ being
a connecting path form $x_0$ to $x_i$. 
We assume that the $\gamma_i$'s have no self-intersections
and no mutual intersection (except at $x_0$).

To each $\gamma_i$, we associate $f_i:=r_{\gamma_i}$.

A \emph{coordinate system} is the (unordered)
$n$-tuple of reflections $\{f_1,\dots,f_n\}$
obtained this way.

We say that a coordinate system is \emph{compatible} with an element
$g\in NC$ if there exists a non-crossing loop $\gamma$ representing $g$,
such that $\Gamma$ is drawn inside $\Int(\gamma)$.
\end{defi}

Coordinate systems are in bijection with isotopy classes of 
planar graphs $\Gamma$ as above (isotopy with fixed vertices).

Saying that $g\in NC$ is compatible with $\{f_1,\dots,f_n\}$
is equivalent to the existence of a permutation $\sigma$ such
that $g=\prod_{i=1}^n f_{\sigma(i)}$. The planar structure
around $x_0$ endows $\{f_1,\dots,f_n\}$ with a natural cyclic
ordering. Once $\{f_1,\dots,f_n\}$ is fixed, choosing a compatible
$g$ is equivalent to the choice of a total ordering refining the 
cyclic ordering (there are $n$ such choices).

Up to isotopy and relabelling of the marked points, 
we may assume that the situation looks like:
$$
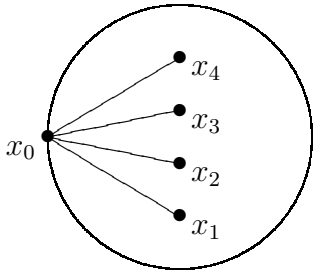
\xy
0;<5pt,0pt>:
(-12,-1) *++{x_0},
(-10,0) *++{\bullet},
(0,-6) *++{\bullet},(0,-2) *++{\bullet},
(0,2) *++{\bullet}, (0,6) *++{\bullet},
(2,-7) *++{x_1}, (2,-3) *++{x_2}, (2,1) *++{x_3}, (2,5) *++{x_4},
(-10,0)="0",
"0";*\xycircle(10,10){},
"0";(0,-6) **\crv{}, 
"0";(0,-2) **\crv{},
"0";(0,6) **\crv{},
"0";(0,2) **\crv{}
\endxy$$
More explicitly, our assumption is that
$x_0=-1$, that the $x_j$ are purely imaginary with
$$-1<\Im(x_1) < \Im(x_2) < \dots < \Im(x_n)<1,$$
and, for each $j$, we consider the affine connecting path
$[x_0,x_j]$ and the associated braid reflection $f_j$.
The coordinate system is then compatible with the element of $NC$
represented by the unit circle.

We have $F_n = \left<f_1,\dots,f_n\right>$.
For any $f\in F_n$, an expression
$f=\prod_{i=1}^m f_{j_i}^{\varepsilon_i}$, with $\varepsilon_i=\pm 1$,
may be obtained as follows. First, find a (possibly self-intersecting)
loop $\gamma$ representing $f$ and drawn inside $D$. Then, following
$\gamma$, write $f_j$ each time it crosses some $[x_j,1]$ moving upwards,
and $f_j^{-1}$ each time it crosses some $[x_j,1]$ moving downwards
(up to perturbation, we may assume that $\gamma$ is transversal to these
segments).

$$
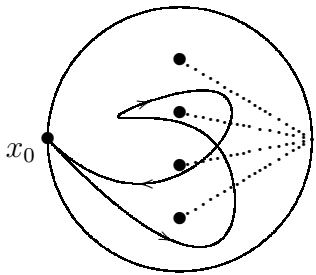
\xy
0;<5pt,0pt>:
(-12,-1) *++{x_0},
(-10,0) *++{\bullet},
(0,-6) *++{\bullet},(0,-2) *++{\bullet},
(0,2) *++{\bullet}, (0,6) *++{\bullet},
(-10,0)="0",(10,0)="1",
"0";*\xycircle(10,10){},
"1";(0,-6) **\crv{~*=<3pt>{.}"1"}, 
"1";(0,-2) **\crv{~*=<3pt>{.}"1"},
"1";(0,6) **\crv{~*=<3pt>{.}"1"},
"1";(0,2) **\crv{~*=<3pt>{.}"1"},
"0";"0" **\crv{(4,-15)&(7,4)&(-10,0)&(5,6)&(4,-2)&(-5,-5)}
?(.06)*\dir{>}?(.5)*\dir{>}?(.9)*\dir{>}
\endxy$$
In the above example, the word is $f_1f_2f_3^{-1}f_2^{-1}$.

A word in the $\{f_1,f_1^{-1},\dots,f_n,f_n^{-1}\}$ is 
\emph{reduced} if the patterns $f_jf_j^{-1}$ and $f_j^{-1}f_j$
never occur.
Any $f\in F_n$ admits a unique expression as a reduced word
in the $\{f_1,f_1^{-1},\dots,f_n,f_n^{-1}\}$.

A loop is \emph{reduced} if the associated word is reduced.
Clearly, any loop in $D$ admits, in its homotopy class,
a reduced loop. More precisely,
this reduced loop may be obtained by a certain ``bigon elimination''
procedure, during which one may avoid introducing self-intersections.
In particular, any non-crossing loop is homotopic to 
a non-crossing reduced loop.

In the next two results, we denote by $g$ the 
(maximal) element of $NC$ represented by $D$.

\begin{lemma}
\label{quadratfrei}
Let $f\in NC_g$. The reduced word associated with $f$
is ``quadratfrei'': it does not contain the patterns
$f_jf_j$ and $f_j^{-1}f_j^{-1}$.
\end{lemma}

\begin{proof}
A picture is worth a thousand words:
$$
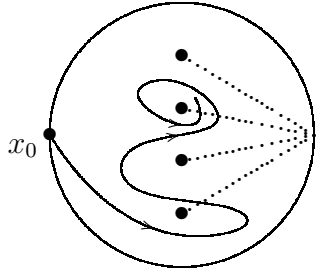
\xy
0;<5pt,0pt>:
(-12,-1) *++{x_0},
(-10,0) *++{\bullet},
(0,-6) *++{\bullet},(0,-2) *++{\bullet},
(0,2) *++{\bullet}, (0,6) *++{\bullet},
(-10,0)="0",(10,0)="1",
"0";*\xycircle(10,10){},
"1";(0,-6) **\crv{~*=<3pt>{.}"1"}, 
"1";(0,-2) **\crv{~*=<3pt>{.}"1"},
"1";(0,6) **\crv{~*=<3pt>{.}"1"},
"1";(0,2) **\crv{~*=<3pt>{.}"1"},
"0";(1,2.8) **\crv{(-5,-8)&(6,-8)&(5,-5)&(-5,-5)&(-5,0)&(5,0)&(0,5)&(-5,3.4)&
(0,0)&(2,1)}
?(.06)*\dir{>}?(.5)*\dir{>}?(.9)*\dir{>}
\endxy$$
\end{proof}

\section{Simple transitivity of Hurwitz actions}

The material in this section is certainly classical, except
the interpretation in terms of Coxeter elements in the universal
Coxeter group.

Choose $g$ a maximal non-crossing element of $F_n$. As we have noted
earlier, it is possible to find a coordinate system $f_1,\dots,f_n$
such that $F_n=\left< f_1,\dots,f_n \right>$ and $g=f_1\dots f_n$.
To fix the notations, we make the standard choice for $g$ and $f_1,\dots,f_n$,
already used in the previous section:
$$
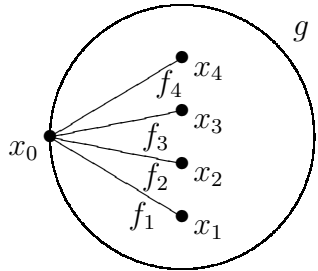
\xy
0;<5pt,0pt>:
(-12,-1) *++{x_0},
(-10,0) *++{\bullet},
(0,-6) *++{\bullet},(0,-2) *++{\bullet},
(0,2) *++{\bullet}, (0,6) *++{\bullet},
(2,-7) *++{x_1}, (2,-3) *++{x_2}, (2,1) *++{x_3}, (2,5) *++{x_4},
(-10,0)="0", (9,8) *++{g}, (-3,-6) *++{f_1}, (-2,-3) *++{f_2},
(-2,0) *++{f_3}, (-1,4) *++{f_4},
"0";*\xycircle(10,10){},
"0";(0,-6) **\crv{}, 
"0";(0,-2) **\crv{},
"0";(0,6) **\crv{},
"0";(0,2) **\crv{}
\endxy$$
Clearly, any expression of $g$ as a product of elements of $R$ must
be of length $n$ (consider the largest abelian quotient of $F_n$).

Thus
$$(f_1,\dots,f_n) \in \Red_R(g).$$

\begin{lemma}
\label{redNC}
Let $\beta \in B_n$, let $(r_1,\dots,r_n) := \beta \cdot (f_1,\dots,f_n)$.
Consider a sequence of integers $j_1,\dots,j_k$ such that
$1\leq j_1 < j_2 < \dots < j_k \leq n$.
Then $r_{j_1}\dots r_{j_k} \in NC_g$.
\end{lemma}

\begin{proof}
The elements $r_1,\dots,r_n$ form a coordinate
system (the $B_n$-action sends coordinate systems to 
coordinate systems). Up to isotopy, all coordinate systems look
the same. This reduces the problem to the case when $\beta=1$, for which
the lemma is obvious.
\end{proof}

\begin{defi}
The universal Coxeter group $W_n$ is defined by the presentation:
$$W_n:= \left< s_1,\dots,s_n | s_i^{2}=1 \right>.$$

We consider the epimorphism $\pi: F_n \twoheadrightarrow W_n, f_j
\mapsto s_j$. We set $T:=\pi(R)$.
Elements of $T$ are called \emph{reflections}.
\end{defi}

We set $c:=\pi(g)$.
It is again easy to see that $$(s_1,\dots,s_n) \in \Red_T(c).$$
The map $\pi^n: \Red_R(g) \rightarrow \Red_T(c)$ is a morphism
of $B_n$-sets (where both sets are equipped with Hurwitz action).

\begin{theo}
\label{orbit}
\begin{itemize}
\item[(1)] The Hurwitz action is simply transitive on $\Red_R(g)$.
\item[(2)] The Hurwitz action is simply transitive on $\Red_T(c)$.
\item[(3)] The map $\pi^n: \Red_R(g) \rightarrow \Red_T(c)$ is
an isomorphism of $B_n$-sets.
\end{itemize}
\end{theo}

The author is grateful to Sang Jin Lee for pointing out that (1)
was already contained in Artin's 1947 article \cite{artin}.

\begin{proof}
The transitivity statement in (1) is Theorem 16 in  \cite{artin}
(although it appears in a formulation closer to ours at the top of
p. 114 of \emph{loc. cit.}). 

Let us prove the transitivity statement in (2) -- our argument is so
similar to Artin's that we could have omitted the proof, but we include
it for the convenience of the reader, who will easily reconstruct the
proof of Artin's Theorem 16.
We start with a remark about normal forms in $W_n$.
This group is a free
product of $n$ cyclic groups of order $2$.
Consider a finite sequence $w:=(a_1,\dots,a_m)$, where each
$a_i$ is taken in $\{s_1,\dots,s_n\}$.
We say that $w$ \emph{represents}
the element $a_1\dots a_m\in W_n$.
We say that $w$ is the \emph{normal form} of $a_1\dots a_m$ if it does
not contain a pattern $s_js_j$  of consecutive equal terms.
When $w$ is a normal form, we say
that $m$ is the \emph{length} of $a_1\dots a_m$.
Clearly, the normal form always exists and is unique. It may actually
be computed with the following non-deterministic procedure.
Start from an arbitrary $w$.
\begin{itemize}
\item[(I)] If $w$ is a normal form, return $w$.
\item[(II)] Otherwise, a least a pattern $s_js_j$ appears. Choose an occurence
and remove the involved terms. Start again with the new (shorter) sequence.
\end{itemize}
A sequence of successive choices in (II) is called an \emph{execution}
of the procedure. Though there are usually
several executions, the end result
is always the (unique) normal form.
The surviving terms in the output come from terms in the input. 
If we choose a particular execution, we say that a given term of $w$
is \emph{untouched} by the execution if it survives it.

Any $t\in T$, being a reflection, may be written
$$(*) \qquad t = u_1 u_{2}
\dots u_{k} s_t
u_{k} \dots u_{2} u_{1}$$
where the $u_i$'s and $s_t$ are
in $\{s_1,\dots,s_n\}$.
We may clearly assume that $(*)$ is a normal form.
We say that $s_t$ is the \emph{content term} of $t$.

Let $(t_1,\dots,t_n)\in \Red_T(c)$.
Considering the largest abelian quotient of $W_n$, one
may observe that, the content terms $s_{t_j}$ satisfy
$\{s_{t_1},\dots,s_{t_n}\}= \{s_1,\dots,s_n\}$.

The normal form of $(t_1,\dots,t_n)$ is $(s_1,\dots,s_n)$.
Let $w$ be the concatenation of the normal forms of $t_1,\dots,t_n$.
Choose an execution of the normal form procedure, applied to $w$.
The output is $(s_1,\dots,s_n)$.
We distinguish two cases:

{\bf Case 1.} The content terms of the normal forms of the $t_j$'s are
untouched by the execution. Write
$$w=(u_1,\dots,u_k, s_1, u_k,\dots,u_1,v_1,\dots,v_l,s_2,
v_l,\dots,v_1,\dots,\dots)$$
(since they are untouched, the content terms must already be in the
order $s_1,\dots,s_n$ in $w$).
The execution rewrites $w$ to $(s_1,\dots,s_n)$ while leaving
$s_1$ untouched. Thus it rewrites $(u_1,\dots,u_k)$ to $()$. Since
$(u_1,\dots,u_k)$ is normal, this implies that $k=0$.
Considering the fragment $v_1,\dots,v_l$ between the unaffected
terms $s_1$ and $s_2$, we conclude that $l=0$, and so on... Thus
$(t_1,\dots,t_n)=(s_1,\dots,s_n)$.

{\bf Case 2.} At least one content term of one the $t_j$ is destructed.
Consider the first iteration of the execution where this happens:
a certain pattern appears, involving (the descendant of)
a content term of at least one
of the $t_j$'s:
denoting by $(a_1,\dots,a_m)$ the word just before this particular 
iteration, we have $a_i=a_{i+1}$ for some $i$,
with $a_i$ or $a_{i+1}$ being the (until then untouched) content term
of one of $t_j$'s. Note that $a_i$ and $a_{i+1}$ may not both
be content terms, because distinct $t_j$'s have distinct contents.
Let us assume that $a_i$ is the content term of some $t_j$.
(The case when $a_{i+1}$ is the content term may be dealt with symmetrically).
Inside $w$, we are interested in the portion involving $t_j$ and $t_{j+1}$:
$$w=(\dots,u_1,\dots,u_k, s, u_k,\dots,u_1,v_1,\dots,v_l,s',
v_l,\dots,v_1,\dots,),$$
where $s$ is the content of $t_j$ and $s'$ the content of $t_{j+1}$.
\begin{lemma}
The length of $s u_k^{-1}\dots u_1^{-1} v_1\dots v_l$ is
$< l-k$ (in particular, $k< l$).
\end{lemma}
\begin{proof}[Proof of the lemma]
From the assumptions, it is easy to see that
the first term $s$ is modified in any execution with input
$(s, u_k,\dots,u_1,v_1,\dots,v_l)$; in particular, this sequence
is not a normal from. Consider an execution with this input.

If $k=0$, we observe that $(v_1,\dots,v_l)$ is a normal form. Since
$(s,v_1,\dots,v_l)$ is not normal, we must have $s=v_1$. The claim holds.

If $k>0$, we observe that both $(s, u_k,\dots,u_1)$
and $(v_1,\dots,v_l)$ are normal forms. 
We must have $u_1$ and $v_1$, and the first step of the execution leads to
$(s, u_k,\dots,u_2,v_2,\dots,v_l)$. We conclude by an easy induction.
\end{proof}

Consider the pair $(t_jt_{j+1}t_j^{-1},t_j)$.
The first reflection is represented by
$$(u_1,\dots,u_k, s, u_k,\dots,u_1,v_1,\dots,v_l,s',
v_l,\dots,v_1,u_1,\dots,u_k, s,
u_k,\dots,u_1).$$
By the lemma, the length of $su_k\dots u_1v_1\dots v_l$ is
$< l-k$. The same property holds for its inverse
$v_l\dots v_1u_1\dots u_ks$.
Thus the length $L$ of $t_jt_{j+1}t_j^{-1}$
satisfies $L <k+(l-k) + 1 + (l-k) + k = 2l+1$.
The total length of $(t_1,\dots,t_{j-1},t_jt_{j+1}t_j^{-1},t_j,t_{j+2},\dots,
t_n)$ is strictly smaller than the total length
of $(t_1,\dots,t_n)$. These two decompositions lie in the same Hurwitz
orbit. One may prove the transitivity part of (2)
by induction on the total length.

The simplicity statement in (1) says that
$$\forall \beta\in B_n, \forall (t_1,\dots,t_n) \in \Red_R(g),
\beta \cdot (t_1,\dots,t_n) = (t_1,\dots,t_n) \Rightarrow \beta=1.$$
Using the transitivity, this statement is equivalent to
$$\forall \beta\in B_n,\beta\cdot (f_1,\dots,f_n) = (f_1,\dots,f_n) \Rightarrow \beta=1,$$
which is nothing but the faithfullness of the standard representation
of $B_n$ in $\Aut(F_n)$, already known to Hurwitz.

Let us now prove the simplicity statement in (2). 
Using transitivity, it is enough to prove that 
$$\forall \beta\in B_n,\beta\cdot (s_1,\dots,s_n)
= (s_1,\dots,s_n) \Rightarrow \beta=1.$$
Let $\beta\in B_n$ such that $\beta\cdot (s_1,\dots,s_n)
= (s_1,\dots,s_n)$.
Let $(r_1,\dots,r_n):=\beta\cdot (f_1,\dots,f_n)$.
Since $\pi^n$ commutes with Hurwitz action,
we have
$\pi^n ((r_1,\dots,r_n) ) = (s_1,\dots,s_n)$,
thus $s_j = \pi(r_j)$ for all $j$. Fix $j\in \{1,\dots,n\}$.
By Lemma
\ref{redNC}, we know that $r_j\in NC_g$.
Consider the normal form $f_{j_1}^{\varepsilon_1} \dots 
f_{j_m}^{\varepsilon_m}$ of $r_j$ in $F_n$. By Lemma
\ref{quadratfrei}, this normal form is ``quadratfrei''.
Thus $s_{j_1} \dots 
s_{j_m}$ is the normal form of $s_j$ in $W_n$. 
Thus $m=1$ and $r_j=f_j^\varepsilon$. Since
$r_j\in R$, we have $\varepsilon =1$.
This holds for any $j$, thus $(r_1,\dots,r_n)=(f_1,\dots,f_n)$.
By (1), we must have $\beta=1$.

(3) follows trivially.
\end{proof}

\begin{coro}
\label{coro1}
There are natural bijections between:
\begin{itemize}
\item[(i)] Maximal strict chains of $NC_g$.
\item[(ii)] Elements of $\Red_R(g)$.
\item[(iii)] Coordinate systems compatible with $g$.
\end{itemize}
More precisely, the map from $(i)$ to $(ii)$ sends a maximal
chain $1= a_0 < a_1< \dots < a_n=g$ to $(a_0^{-1}a_1,\dots,a_{n-1}^{-1}a_n)$,
and the map from $(ii)$ to $(iii)$ send $(t_1,\dots,t_n)$ to
$\{t_1,\dots,t_n\}$.
\end{coro}

\begin{proof}
Consider the classical interpretation of 
$B_n$ as the mapping class group of the $n$-punctured
disk, fixing the outer circle.

By Lemma \ref{stupidlemma}, maximal strict chains of $NC_g$ are
represented by chains of concentric non-crossing loops in $D$, of
stricly increasing height.
Isotopy classes of such data clearly form a single $B_n$-orbit.

Similarly, coordinates systems drawn inside $D$ form
a single $B_n$-orbit.

The corollary then follows from the fact that $\Red_R(g)$ is a single
Hurwitz orbit, and that the natural maps with the above objects are
$B_n$-equivariant.
\end{proof}

\begin{coro}
\label{corocon}
Denote by $R_g$ the subset of $R$ consisting of elements
which may appear in some sequence in $\Red_R(g)$.
Denote by $T_c$ the subset of $T$ consisting of elements
which may appear in some sequence in $\Red_T(c)$.
Then $R_g = R \cap NC_g$. Moreover, $\pi$ induces a
bijection $R_g \simeq T_c$.
\end{coro}

Note that $\pi$ does not induce a bijection from $R$ to $T$.
Also, the injectivity of $R_g \simeq T_c$
is \emph{a priori} stronger than the injectivity of $\pi^n:\Red_R(g)
\rightarrow \Red_T(c)$ from the theorem.

\begin{proof}
The statement $R_g = R \cap NC_g$ is already in Lemma \ref{redNC}.
Using the theorem, we note that $T_c = \pi(R_g)$.
We are left with having to prove the injectivity.
First, we observe that the fiber of $R_g \rightarrow T_c$ over
$s_1$ is a singleton (it follows from Lemma \ref{quadratfrei}).
By transitivity, all fibers have the same cardinal.
\end{proof}

\section{Quasi-Garside structure}

\begin{defi}
We denote by $F_n^+$ the submonoid of $F_n$ generated by $R$.
We endow $F_n^+$ with the divisibility partial ordering:
for all $f,g\in F_n^+$,
$f\preccurlyeq g \stackrel{\text{def}}{\Longleftrightarrow}
\exists h\in F_n^+, fh=g$.
\end{defi}

Note that, since $R$ is a an union of conjugacy classes, 
$\exists h\in F_n^+, fh=g \Leftrightarrow \exists h\in F_n^+, hf=g$.
We do not have to distinguish left divisibility from right divisibility.

\begin{lemma}
The restriction of $\preccurlyeq$ to $NC$ coincides with $\subseteq$.
\end{lemma}

\begin{proof}
Let $f,g\in NC$.

It is constructively clear that $f \subseteq g$ implies $f \preccurlyeq g$.

Conversely, if $f \preccurlyeq g$, then a reduced $R$-decomposition
$(r_1,\dots,r_k)$ of $f$ may be extended to a reduced $R$-decomposition
$(r_1,\dots,r_l)$ of $g$. By Lemma \ref{redNC},
$r_1\dots r_k \in NC_g$.
\end{proof}

In \cite[Definition 0.5.1]{dualmonoid}, a
Garside monoid was defined as a monoid
$M$ satisfying a certain number of axioms; one of these axioms concerns
the existence of a ``balanced'' element $\Delta\in M$ whose set
of left/right divisors is finite and generates $M$.

For many applications, one may work in a slightly generalised context:
by \emph{quasi-Garside} monoid, we mean a monoid satisfying all
axioms of \cite[0.5.1]{dualmonoid}, except that we do not require the
set of divisors of $\Delta$ to be finite.

\begin{theo}
Let $g$ be a maximal element of $NC$. 
Let $M_g$ be the submonoid of $F_n$ generated by
$\{r\in R | r \preccurlyeq g\}$.
Then $M_g$ is a quasi-Garside monoid with Garside element $g$ and
set of simples $NC_g$. 
\end{theo}

\begin{proof}
Set $P_g:=\{r\in R | r \preccurlyeq g\}$.
Using a straightforward analog of \cite[Theorem 0.5.2]{dualmonoid}
where the finiteness condition is removed, we only have to prove
that $(P_g,\preccurlyeq)$ is a lattice.

Using the last lemma, we see that any element of $NC_g$ lies
in $P_g$; conversely, using Corollary \ref{coro1}, we see
that any element of $P_g$ belongs to $NC_g$; using again the last
lemma, we have $(P_g,\preccurlyeq)=(NC_g,\subseteq)$. By Theorem 
\ref{theolattice}, the latter is a lattice.
\end{proof}

The free group being easy enough to study with the classical point
of view (with its presentation with $n$ generators and $0$ relations)
that what brings the above quasi-Garside structure may seem futile:
for example, 
we have a new presentation with an infinity of generators (reflections
in $NC_g$) and
an infinity of relations of length $2$ (the relations $rr'=r''r$,
whenever $r,r'\in NC_g$ satisfy $rr'\in NC_g$ and $r''=rr'r^{-1}$),
with a solution to the
word and conjugacy problem...
The main interest of this quasi-Garside structure is that it
fits in a general pattern, formalised in the conjectures below,
and also that it is useful to understand geometric
aspects of complex reflection groups, as it will appear in the
sequel \cite{combi}.

\section{Conjectures}

As announced in the introduction, our conjectures apply to two different
settings:
\begin{itemize}
\item[(a)] either $(W,S)$ is a Coxeter system; we assume that $n:=|S|$ is
finite (but $W$ may be 
infinite); we denote by $T$ the set of reflections
in $W$ (arbitrary conjugates in $W$ of elements of $S$);
we consider the associated Artin group $B:=A(W,S)$ (we will use
bold fonts to refer to the formal copy of $S$ generating $B$);
we denote by $R$ the set of ``braid reflections'' (arbitrary conjugates
in $B$ of elements of $\BS$);
\item[(b)] or $W$ is an irreducible complex reflection group of rank $n$
generated by involutive reflections;
we assume that it is ``well-generated'', \ie, it may be generated by $n$
reflections; we denote by $T$ the set of all reflections in $W$;
we consider the generalised braid group $B:=B(W)$, defined in \cite{bmr}
as the fundamental group of the space of regular orbits; we denote
by $R$ the set of ``braid reflections'' (``generators-of-the-monodromy'')
in $B$.
\end{itemize}

In both settings, there is a natural map $p: B \twoheadrightarrow W$.

\begin{defi}
A \emph{Coxeter element} is, depending on the setting:
\begin{itemize}
\item[(a)] the conjugate in $W$ of a product $s_1\dots s_n$,
for a certain numbering $S=\{s_1,\dots,s_n\}$;
\item[(b)] an element $c\in W$ such that
$\ker(c-e^{\frac{2i\pi}{d_n}})\neq 0$, where
$d_n$ is the largest invariant degree of $W$.
\end{itemize}
A \emph{braid Coxeter element} is, depending on the setting:
\begin{itemize}
\item[(a)] the conjugate in $B$ of a product $\bs_1\dots \bs_n$,
for a certain numbering $\BS=\{\bs_1,\dots,\bs_n\}$;
\item[(b)] an element $g\in B$ such that $g^{d_n}=\pi$, where
$d_n$ is the largest invariant degree of $W$, and $\pi$ is the standard
``full-turn'' element in the center of $B$ (\cite{bmr}).
\end{itemize}
\end{defi}

Clearly, in the situation (a), $p$ maps braid Coxeter elements to Coxeter
elements. This also holds in $(b)$ (\cite{combi}).

An important issue is that, in situation (a), there are usually
several conjugacy classes of (braid) Coxeter element. 
However, when the Coxeter graph is a tree, there is a unique
conjugacy class (\cite{LIE}, p. 117).
In the situation of the free group, there are many conjugacy classes,
but they are group-theoretically undistinguishable, since the full symmetric
group acts by diagram automorphisms.

In our conjectures, only the conjugacy class of the braid Coxeter element
matters.

\begin{conj}
\label{conj1}
There exists a braid Coxeter element $g\in B$
such that, setting $c:=p(g)$, we have:
\begin{itemize}
\item[(1)] The Hurwitz action is transitive on $\Red_R(g)$.
\item[(2)] The Hurwitz action is transitive on $\Red_T(c)$.
\item[(3)] The map $p^n$ induces an isomorphism of $B_n$-sets
from $\Red_R(g)$ to $\Red_T(c)$.
\item[(4)] The map $p$ induces a bijection from the set $R_g$ of
reflections appearing in $\Red_R(g)$ to the set $T_c$ of braid reflections
appearing in $\Red_T(c)$.
\end{itemize}
\end{conj}

In the case of the universal Coxeter group $W_n$ and its braid group
$F_n$, the conjecture is proved above (Theorem \ref{orbit} and
Corollary \ref{corocon}). 
That the action is then \emph{simply} transitive and not just
transitive is specific to this
case.

When $W$ is a finite Coxeter group,
most of the conjecture
is proved in \cite{dualmonoid}: (2) is \emph{loc. cit.} Proposition 1.6.1, and
a weaker form of (3) and (4) are consequences of Fact 2.2.4; however, 
no description of $\Red_R(g)$ is given (only a specific $B_n$-orbit
is considered, it is not proved to be the full $\Red_R(g)$).

When $W$ is the Coxeter group of type $\widetilde{A_n}$,
this follows from \cite[Proposition 3.4]{digne}. Note that Digne proves a more
general result: the transitivity is true for all braid Coxeter
elements. The above conjecture
is certainly not optimal (see for example Digne's Conjecture 1.1).
Actually, in view of \cite[Theorem 3.16]{brieskorn}
(and the discussion following this result on p. 87), it is
tempting to formulate a more general conjecture, not only
applying to Coxeter elements but to elements whose reduced
decompositions involve generating sets. However, since
we have neither interesting examples nor applications,
we stay with the above conjecture, which interests us
in connection with our second conjecture below.

Given any braid Coxeter element $g$,
consider the positive presentation with 
set of generators $R_g$ and relations $rr'=r''r$ whenever there exists
an element of $\Red_R(g)$ starting by $(r,r',\dots)$ and $r''=rr'r^{-1}$.
Let $M_g$ be the monoid defined by this presentation; let $B_g$ be
the group defined by this presentation.

Since the relations $rr'=r''r$ hold in $B$, $B$ is \emph{a priori} a quotient
of $B_g$. In setting (a), it is easy to see that the 
defining relations of $B$ are consequences of the Hurwitz relations,
thus that $B_g\simeq B$. One may prove the
similar statement in setting (b) (\cite{combi}).

Points (1) and (2) of the above conjecture express that $M_g$ coincides
with the
monoids associated with the triples $(B,R,g)$ and $(W,T,c)$, as in
\cite[Section 0.4]{dualmonoid}. With the obvious analog
of \cite[Theorem 0.5.2]{dualmonoid}, the next conjecture is
the key ingredient to prove that $M_g$ is a quasi-Garside monoid.

\begin{conj}
\label{conj2}
Denote by $B_+$ the submonoid of $B$ generated by $R$.
Denote by $\preccurlyeq$ the relation on $B_+$ defined by
$b\preccurlyeq b'$ if and only if $b^{-1}b'\in B_+$.
For any $b\in B_+$, set $P_b :=\{b'\in B_+| b'\preccurlyeq b\}$.
There exists a braid Coxeter element $g\in B$ satisfying Conjecture 
\ref{conj1} and such that $(P_g,\preccurlyeq)$ is a lattice.
\end{conj}

Again, this is known for spherical Artin types, \cite{dualmonoid},
and affine type $\tilde{A}$, \cite{digne}, and in $F_n$ as it was proved above.
The most mysterious aspect
is that the lattice does \emph{not} hold for all Coxeter elements:
indeed, Digne's striking Proposition 5.5 shows that, in $\widetilde{A_{n-1}}$,
it holds only when the braid Coxeter element is a product of the 
generators according to the cyclic order on the diagram.
We have no good hint on how to characterise suitable braid Coxeter
elements in setting (a). In setting (b), all choices are conjugate.

Among possible applications, we observe that
braid groups satisfying conjectures \ref{conj1} and \ref{conj2}
have cohomological dimension smaller or equal to $n$,
since the construction of \cite{charney} of a simplicial $K(\pi,1)$
for Garside groups clearly extends to quasi-Garside groups
(the obtained $K(\pi,1)$ still being of dimension $n$, but no longer
necessarily finite).


\begin{thebibliography}{DDRW}      

\bibitem[A]{artin} E. Artin,
{\em Theory of braids}, Ann. of Math. (2) {\bf 48} (1947), 101-126.

\bibitem[B1]{dualmonoid} D. Bessis,
{\em The dual braid monoid},
Ann. Sci. \'Ecole Norm. Sup. {\bf 36} (2003), 647--683. 

\bibitem[B2]{combi} D. Bessis,
{\em On well-generated complex reflection groups},
in preparation.

\bibitem[Br]{brieskorn} E. Brieskorn,
{\em Automorphic sets and braids and singularities}, in 
{\em Braids}, Contemporary Mathematics {\bf 78}, American Mathematical
Society, 1988, 45--117.

\bibitem[BC]{BC} D. Bessis, R. Corran,
{\em Garside structure for the braid group of type $(e,e,n$)},
in preparation.

\bibitem[BKL]{BKL} J.~Birman, K.~H.~Ko, S.~J.~Lee,
{\em A new approach to the word and conjugacy problem in the braid
groups}, Adv. Math. {\bf 139} (1998), no.2, 322--353.


\bibitem[BMR]{bmr} M.~Brou\'e, G.~Malle et R.~Rouquier,
        {\em Complex reflection groups, braid groups, Hecke algebras},
        J. reine angew. Math. {\bf 500} (1998), 127--190. 


\bibitem[CMW]{charney} R. Charney, J. Meier, K. Whittlesey,
{\em Bestvina's normal form complex and the homology of
Garside groups}, preprint (2001).

\bibitem[D]{digne} F. Digne,
{\em Pr\'esentations duales des groupes de tresses de type affine
$\tilde{A}$}, preprint.


\bibitem[DDRW]{panobraid} P.~Dehornoy, I.~Dynnikov, D.~Rolfsen, B.~Wiest,
{\em Why are braids orderable?},
Panoramas et Synth\`eses {\bf 14} (2002),
Soci\'et\'e Math\'ematiques de France, Paris.

\bibitem[DP]{dehpa} P.~Dehornoy, L.~Paris, {\em Gaussian groups and Garside
groups, two generalizations of Artin groups},
Proc. of London Math. Soc. {\bf 79} (1999), 569--604.

\bibitem[LIE]{LIE} N. Bourbaki,
{\em Groupes et alg\`ebres de Lie}, chapitres IV, V et VI,
Hermann, 1968.

\end{thebibliography}
\end{document}